\documentclass{article}
\usepackage[utf8]{inputenc}

\usepackage[margin = 1in]{geometry}
    \usepackage[numbers]{natbib}
	\usepackage[toc,page]{appendix}
    \usepackage[nottoc]{tocbibind}
    \usepackage[english]{babel}
    \usepackage{marvosym}
	\usepackage{amssymb, amsmath,bm,graphicx}
	\usepackage{setspace}
	\usepackage{pdfpages}
	\usepackage{mathrsfs,amsthm}
    \usepackage{authblk}
	\usepackage[colorlinks,citecolor=blue,urlcolor=blue,linkcolor = blue]{hyperref}
	\usepackage[normalem]{ulem}
	\usepackage{subcaption}
	\usepackage{comment}
	\usepackage{multirow}
	\usepackage{anyfontsize}
	\usepackage{graphicx,psfrag,epsf}
	\usepackage{float}
	\usepackage{thmtools}
	\usepackage[labelfont=bf]{caption}
	\usepackage{tabularx}
	\usepackage{url}
	\usepackage{booktabs}
	\usepackage{hyperref}
	\usepackage{color}
	\usepackage{xcolor}
	\usepackage{caption}
	\usepackage{bm,bbm,bigstrut,geometry,pdflscape}
	\usepackage[mathscr]{euscript}
	\usepackage{bigstrut,enumerate}
	\usepackage[utf8]{inputenc}
	\usepackage{enumitem}
	\usepackage[toc,page]{appendix}
    \usepackage[noabbrev,capitalize]{cleveref}
	\usepackage[ruled,vlined]{algorithm2e}

    \theoremstyle{plain}
	\newtheorem{theorem}{Theorem}[section]
	\newtheorem{lemma}[theorem]{Lemma}
	\newtheorem{remark}[theorem]{Remark}
	
	\newtheorem{assumption}{Assumption}
	\newtheorem{defn}[theorem]{Definition}

\numberwithin{equation}{section}
\numberwithin{theorem}{section}
\numberwithin{lemma}{section}
\numberwithin{corollary}{section}
\numberwithin{remark}{section}

    \usepackage{color,soul}

\newcommand{\R}{\mathbb{R}}
\newdimen\slantmathcorr
\def\oversl#1{
\setbox0=\hbox{$#1$}
\slantmathcorr=\wd0
\hskip 0.2\slantmathcorr \overline{\hbox to 0.8\wd0{%
\vphantom{\hbox{$#1$}}}}
\hskip-\wd0\hbox{$#1$}
}

\def\undersl#1{
\setbox0=\hbox{$#1$}
\slantmathcorr=\wd0
\underline{\hbox to 0.8\wd0{%
\vphantom{\hbox{$#1$}}}}
\hskip-0.8\wd0\hbox{$#1$}
}

\newcommand{\bP}{\mathbb{P}}

\usepackage{soul}

\usepackage{bm,bbm}
\usepackage{dsfont}

\renewcommand{\epsilon}{\varepsilon}

\newcommand{\E}{\mathbb{E}}

\usepackage{xcolor}

\SetKwInOut{Parameter}{Parameter}
\usepackage{authblk}

\usepackage{commath}
\newcommand{\p}[2][{}]{\mathbb{P}_{#1} \left({#2} \right)}
\newcommand{\e}[2][{}]{\mathbb{E}_{#1} \sbr{#2}}
\newcommand{\expp}[1]{\exp \left( {#1} \right)}
\newcommand{\paranth}[1]{\left( {#1} \right)}
\newcommand{\N}{\mathcal{N}}

\newcommand{\define}{\triangleq}
\definecolor{Milad}{RGB}{200, 0, 100}

\title{Exponential tail bounds and Large Deviation Principle for Heavy-Tailed U-Statistics\footnote{This work was partially supported by the National Science Foundation under Grant 1811614. }}
\author{Milad Bakhshizadeh}
\affil{Stanford University}

\date{}

\begin{document}

\maketitle

\begin{abstract}
    We study deviation of U-statistics when samples have heavy-tailed distribution so the kernel of the U-statistic does not have bounded exponential moments at any positive point.  We obtain an exponential upper bound for the tail of the U-statistics which clearly denotes two regions of tail decay, the first is a Gaussian decay and the second behaves like the tail of the kernel.  For several common U-statistics, we also show the upper bound has the right rate of decay as well as sharp constants by obtaining rough logarithmic limits which in turn can be used to develop LDP for U-statistics.  In spite of usual LDP results in the literature, processes we consider in this work have LDP speed slower than their sample size $n$.
\end{abstract}

\section{Introduction}
Suppose $X_1,X_2,\ldots ,X_n \overset{iid}{\sim} \mu$ where $\mu$ is a probability measure supported on $\R$. We consider the following  U-statistic of order $m$ with a kernel function $h(\cdot):\R^m\to\R$ which is symmetric about its arguments:
$$U_n \define \frac{1}{{n\choose m}}\sum_{1\leq i_1<\ldots <i_m\leq n} h(X_{i_1},\ldots ,X_{i_m}).$$

We study the decay of $\p{\abs{U_n - \e{U_n}} > t}$ as a function of $t, n$, and $m$, in a setup which is rarely addressed in the literature and is characterized by a couple of key assumptions: heavy-tailed distribution for $h(X_1, ..., X_m)$, and large values for $t$.  We explain the role of each assumption in more details below.

When $h(X_1, ..., X_m) \define h$ has a heavy-tailed distribution, its Moment Generating Function (MGF) is not bounded at any positive point.  This breaks many concentration results and violates Cramer's condition which is required for common results that determine large deviation behavior. While deviation of U-statistics has been studied extensively in the light-tailed regime \cite{arcones1992large, eichelsbacher1995large, eichelsbacher2005refined, korolyuk2013theory, lee2019u, peaucelle2004efficiency, wang1994large}, the same for heavy-tailed U-statistics is relatively under-explored. In this paper, we aim to focus on the heavy-tailed regime.

Moreover, when kernel $h$ is assumed to have a heavy tail, for fixed values of $n$ and $m$,  $\p{\abs{U_n - \e{U_n}} > t}$ has two different behaviors: 1) for small values of $t$ it decays like a Gaussian tail, and 2) for $t$ larger than the zone of normal convergence, i.e. large deviation zone, it has a decay slower than normal distribution.  In spite of the first region which has been studied in the literature largely, see \citep{lee2019u, korolyuk2013theory} and the references therein, there are little documented information about the tail behavior in the second region.  Several setups have been developed to study the behavior of the tail. We mention some of them below to discuss the setup that results of the current paper belong to.
\begin{enumerate}
    \item Finite sample concentration inequalities that give upper bounds for $\p{\abs{U_n - \e{U_n}} > t}$ for all values of n.
    \item Asymptotic distribution studies that find suitable scaling $a_n$ and limiting distribution $D$ for which $a_n \abs{U_n - \e{U_n}} \xrightarrow{d} D$.  $a_n = c \sqrt{n}$, and $D = \mathcal{N}(0, 1)$ is the well-known CLT that holds for non-degenerate U-statistics.
    \item Berry–Esseen type inequalities that seek for uniform upper bound $\abs{\p{a_n \abs{U_n - \e{U_n}} > t} - \phi(t)}$ for all $t \in \mathcal{C}$, where $\mathcal{C} \subseteq \mathbb{R}$.
    \item Large deviation studies that look for convergence speed $b_n$ and rate function $f(t)$ for which we have $\lim\limits_{n \to \infty} \frac{- \log \p{\abs{U_n - \e{U_n}} > t}}{b_n} = f(t), \; f(t) > 0$ for large $t$.
\end{enumerate}

In this work, we try to shed some light on the undocumented facts about the deviation of U-statistics in setups 1 and 4 listed above, the concentration inequality, and the large deviation setups.  In particular, we develop a general concentration inequality that holds for all U-statistics with finite second moments.  Moreover, we characterize a class of kernels for which that concentration inequality is sharp enough to yield the large deviation limit.  The byproduct of such sharpness analysis is obtaining exact convergence speed $b_n$ and a closed form for the rate function in terms of the tail of the kernel $h$. For heavy-tailed kernels, which are considered in this work, we usually obtain $b_n \ll n$, e.g $b_n = n^\alpha, \; \alpha < 1$.  This distinguishes our results from several previous works that tried to determine tail behavior in the asymptotic Gaussian regime, i.e. $b_n = n$.

\subsection{Related works}

Large deviations of U-statistics have been studied in several previous works under a variety of assumptions.
\citet{hoeffiding1948class} formally introduced U-statistics and showed with the right scaling they converge to a Normal distribution.  Follow up studies offered bounds for the deviation of U-statistics from their Gaussian limits (see \cite{lee2019u, korolyuk2013theory} and the references therein).

\citet{gine2000exponential} offers upper bound for moments of U-statistics.  There is a relatively straightforward connection between moment inequalities and exponential bounds for the tail (see \cite{vershynin2018high}).  Utilizing such connection, \cite{gine2000exponential} obtains exponential tail bounds.  However, these bounds do not show Gaussian decay when the deviation amount $t$ is within the Gaussian zone.  Therefore, the authors also offer improvements for their inequalities in the case of completely degenerate kernels of order $m = 2$.

Recently, \citet{chakrabortty2018tail} considered degenerate  kernels of order 2 and samples from heavy-tailed subWeibull distributions and obtained exponential tail bounds for such U-statistics.
However,
the specific form they used for the U-statistics seems 
restrictive.
They define 
$$U_n = \sum\limits_{i \neq j} \phi(Y_i) w_{i, j}(X_i, X_j) \psi(Y_j),$$
and impose uniform boundedness assumption on $w_{i, j}$.  Hence, the kernel can be unbounded only through product function $\phi(Y_i) \psi(Y_j)$.  It is not clear to us if results of \cite{chakrabortty2018tail} can recover exponential bounds for general non-degenerate kernels like ones named in Lemma \ref{lem: kernels with ldp} of the current work.  

One property of U-statitics of order $m \geq 2$ that makes them more challenging than the case $m=1$, i.e. iid sums, is containing dependent terms in the sum. There are a couple of ways to handle such dependencies, decoupling and martingale inequalities.  Each direction points to a path that is relatively different from another.  Both directions have been explored to develop tail bounds for U-statistics. 

\citet{de1992decoupling} provides a decoupling tool that helps one to get rid of dependencies of summands in the U-statistics.  This tool has been utilized by several later works to obtain exponential concentration bounds, e.g. \cite{liu2020large, bretagnolle1999new, chakrabortty2018tail}. While the decoupling technique is proved very powerful in handling dependencies, there is an extra $8$ factor on the upper bound it offers which usually makes the exponential bound it provides off by a constant factor in its exponent.

Another approach to obtain exponential inequalities for U-statistics is to leverage inequalities for (sub/super) martingales such as those in \cite{fan2012large}. \citet{houdre2003exponential} and some of their references pursue this direction.  Using \citet{de1992decoupling} decoupling and \citet{talagrand1996new} inequality is common in such approach.  Nevertheless, in this approach, the statements get algebraically complicated soon and several constants given by $\sup$ or expectation of partial sums get involved. As a result, many bounds obtained by martingale inequalities are restricted to order 2 or degenerate kernels.  Moreover, constants are not usually sharp for the similar reason discussed in utilizing the decoupling technique above.  

In addition to tail bounds, several works have been devoted to showing U-statistics satisfy Large Deviation Principal (LDP).  It is common for works in this direction to assume the kernel has finite exponential moments, hence they exclude distributions we are focusing on in the current manuscript.

\citet{dasgupta1984large} tried to show large deviation behavior of U-statistics is determined by their {conditional expectation given one summand.} At the first glance this sounds a reasonable claim since the asymptotic normality of U-statistics is tightly related to such conditional expectations \cite{hoeffiding1948class}.  However, this claim has been disapproved later \cite{baringhaus2002large, wang1994large}.

\citet{nikitin2001rough} study large deviation of U-statistics with bounded non-degenerate kernels.  It also claims no large deviation results have been known by the time the work was published.
\citet{arcones1992large} developes a decomposition technique for the case of a 2-variables kernel with finite MGF that enables one to prove LDP.  Nevertheless, extension of this technique to higher order kernels is not straightforward. 
\citet{eichelsbacher1995large} consider both U-statistics and V-statistics, which includes terms with repeated samples in the sum, given by a kernel function with bounded MGF. In this case, they show satisfying LDP for U-statistics and the corresponding V-statistics is equivalent with the same rate function.  Lemma 3.1 from \cite{eichelsbacher1995large}, which allows one to bound MGF of a U-statistic, is one the most important tools we used in the current work.
\citet{eichelsbacher2005refined} provides sufficient conditions under which a U-statistic would satisfy an LDP under weak Cramer condition.  It utilizes a representation of a bi-variable kernel as infinite sums of separable functions in order to reduce the problem of the large deviation for U-statistics to the one for sums of iid variables.  Then, it applies contraction theorem to prove LDP for the U-statistics.  

The current work focuses on heavy-tailed distributions which are mostly excluded from large deviation studies above.  The exponential tail bound given in Section \ref{sec: upper bound} includes all kernels with finite variances regardless of their order, boundedness,  or degeneracy.  The upper bound we offer is simple enough to reveal both Gaussian and heavy tail decay regions immediately.  All parameters in the statement of the upper bound have known limits when $n \to \infty$.  Moreover, we offer ready-to-use tools to handle them for fixed sample size $n$ in Section \ref{sec: paramters of upper bound}.  In addition, to the best of our knowledge, the rough logarithmic limits, with speeds slower than $n$, obtained in Section \ref{sec: LDP} have not been documented in earlier works.

\subsection{Some Applications} \label{sec: application}

U-statistics appear in several statistical analyses including point estimation, such as population variance \cite{hoeffiding1948class},  hypothesis testing \cite{nikitin1995asymptotic}, and statistical mechanics \cite{kiessling2012onsager}.  As a self-averaging function, it is quite natural to ask how fast a U-statistic concentrates around its mean.

In particular, for the asymptotic analysis of hypothesis testing, the speed and the rate function for the large deviation of such statistics are  important. For instance, obtaining the values of Bahadur efficiency and Bahadur exact slope requires having LDP for the test function \cite{nikitin1995asymptotic}.  
Moreover, the large deviation principle and exponential concentrations for U-statistics play substantial roles in the study of macroscopic limits of interacting particle systems and in general statistical mechanics \cite{liu2020large}.
Our result enables one to extend such theories when samples are drawn from a heavy-tailed distribution.

\section{Main results}
The main result of this paper is two-fold. First, in Section \ref{sec: upper bound} we develop a general upper bound for the tail probability of U-statistics whose kernels have bounded variances.  Second, in Section \ref{sec: LDP}, we show for several kernel functions this upper bound is tight enough to determine the large deviation behavior of the U-statistics.  If a centered kernel $h:\mathbb{R}^m \to \mathbb{R}$ satisfies the criteria discussed in Section \ref{sec: LDP} and $U_n$ denotes the U-statistic corresponding to $h$ with $n$ samples, we show in Lemma \ref{lem: ldp} that $\p{U_n > t}$ and $\p{h > \frac{n}{m} t}$ have the same asymptotic behavior in logarithmic sense, i.e.
\begin{equation*}
    \lim_{n \to \infty} \frac{\log \p{U_n > t}}{\log \p{h > \frac{n}{m} t}} = 1.
\end{equation*}
Therefore, both the convergence speed and the rate function for large deviation of $U_n$ are determined by the tail of kernel $h$.  

In addition, the lower bound developed in Section \ref{sec: lowerbound} reveals the event that is responsible for large deviations of $U_n$.  Indeed, one of the $X_1, ..., X_n$ gets large enough to make all the summands it is contributing to exceed $\frac{n}{m} t$.  There will be ${n - 1} \choose {m - 1}$ such terms in the expansion of $U_n$.

\subsection{Finite sample upper bound} \label{sec: upper bound}

Without loss of generality, we will operate under the assumption that the kernel is centered.

\begin{assumption} \label{assump: centered kenrnel}
Suppose $h$ is centered, i.e.
    \begin{equation*}
        \E h(X_1,\ldots ,X_m)=0.
    \end{equation*}
\end{assumption}

We need to define rate functions $I, J$ as below.
\begin{defn} \label{def: I, J}
Let $I, J: \mathbb{R}^{\geq 0} \to \mathbb{R}^{\geq 0}$  be two non-decreasing functions that upperbound right tails of $h$ and $\abs{X_i}$ exponentially.  In other words:
\begin{align} \label{eq: I def}
    \p{h(X_1, ..., X_m) > t} \leq \expp{-I(t)}, \quad \forall t \geq 0,
    \\ \label{eq: J def}
    \p{\abs{X_i} > t} \leq \expp{-J(t)}, \quad \forall t \geq 0.
\end{align}
\end{defn}
Note that one can simply take $I(t) = - \log \p{h(X_1, ... , X_m) > t}, \; J(t) = - \log \p{\abs{X_i} > t}$.  The whole point of Definition \ref{def: I, J} is to allow one to work with possibly simpler upper bounds.

Below Lemma is a simpler modified version of Lemma 3.1 from \cite{eichelsbacher1995large}.
\begin{lemma}\label{lem:mgfbd}
Let $k \define \lfloor \frac{n}{m}\rfloor$. Then, given any $L\geq 0$ and $\lambda\geq 0$, the following holds:
$$\E\left[\exp\left(\lambda \frac{1}{{n\choose m}}\sum_{1\leq i_1<\ldots <i_m \leq n} h_L(X_{i_1},\ldots ,X_{i_m}) \right)\right]\leq 
\e{ \expp{\frac{\lambda}{k} h_L(X_1, ..., X_m)}}^k,
$$
where 
$h_L(X_{i_1},\ldots ,X_{i_m}) \define h(X_{i_1},\ldots ,X_{i_m})\mathbf{1}(h(X_{i_1},\ldots ,X_{i_m})\leq L)$.
\end{lemma}

\begin{proof}
Define $B(X_1, ..., X_n) = \frac{1}{k} (h_L(X_1, ..., X_m) + h_L(X_{m + 1}, ..., X_{2 m}) + ... + h_L(X_{km - m + 1}, ..., X_{k m}))$, then we have
\begin{equation*}
    \frac{1}{{n \choose m}} \sum_{i_1 < ... < i_m} h_L(X_{i_1}, ..., X_{i_m}) = \frac{1}{n !} \sum_{\sigma \in S_n} B(X_{\sigma(1)}, ..., X_{\sigma(n)}).
\end{equation*}
Since $h_L \leq L$ is bounded, its moment generating function is finite at any positive point, hence we obtain

\begin{align*}
    \E\left[\exp\left(\lambda \frac{1}{{n\choose m}}\sum_{1\leq i_1<\ldots <i_m \leq n} h_L(X_{i_1},\ldots ,X_{i_m}) \right)\right]
    & = 
    \E\left[\exp\left( \frac{\lambda}{{n!}}\sum_{\sigma} B(X_{\sigma(1)}, ..., X_{\sigma(n)}) \right)\right]
    \\ & \overset{*}{\leq} 
    \frac{1}{n !} \sum_\sigma \e{\expp{\lambda B(X_{\sigma(1)}, ..., X_{\sigma(n)})}}
    \\ & = 
    \e{ \expp{\lambda B(X_1, ..., X_n)} }
    \\ & =
    \e{ \expp{\frac{\lambda}{k} h_L(X_1, ..., X_m)}}^k.
\end{align*}

To obtain inequality marked by $*$, we used the fact $\expp{\frac{1}{n!} \sum \lambda B \circ \sigma} \leq \frac{1}{n!} \sum \expp{B \circ \sigma}$ by convexity of the exponential function.

\end{proof}

For the sake of simplicity we drop arguments of $h(X_1, ..., X_m), h_L(X_1, ..., X_m)$ and only write $h, h_L$ 

\begin{lemma}[Lemma 1 of \cite{Bakhshizadeh2020}] \label{lem:hmgfbd}
If $\e{h} = 0$, for any $\eta, L \geq 0$ we have
\begin{equation*}
     \e{\exp (\eta h_L)} \leq \exp \left( \frac{v(L, \eta)}{2} \eta^2 \right),
\end{equation*}
where
$v(L, \eta) \define  \e{h_L^2 \mathbf{1}(h \leq 0) + h_L^2 \exp(\eta h_L) \mathbf{1}(h > 0)}$.
\end{lemma}

\begin{theorem} \label{thm:general upper bound}
Under Assumption \ref{assump: centered kenrnel}, for any $0 < \beta \leq 1, \; t \geq 0$
\begin{equation} \label{eq:upper bound}
 \p{U_n > t} \leq \expp{- \frac{k t^2}{2 v(kt, \beta \frac{I(kt)}{kt})}} + \expp{- \beta I(kt) \max(\frac{1}{2}, c(t, \beta, k))} + {n \choose m} \expp{-I(kt)},   
\end{equation}
where 
$v(\cdot, \cdot)$ is the same as in Lemma \ref{lem:hmgfbd}, 
$c(t, \beta, k) \define 1 - \frac{\beta}{2 t} \frac{I(kt)}{kt} v(kt, \beta \frac{I(kt)}{kt})$, and $ k = \lfloor \frac{n}{m} \rfloor$.
\end{theorem}

\begin{proof}
Let show the U-statistic with kernel $h$ by   $U_n(h)$
\begin{alignat*}{2}
\p{U_n (h) > t} & \leq \p{U_n(h_L) > t} + \p{\exists i_1, ..., i_m, \quad h(X_{i_1}, ..., X_{i_m}) > L }    
\\ & \leq
\expp{- \lambda t} \e{\expp{\lambda U_n(h_L)}} + {n \choose m} \expp{- I(L)}
\\ & {\leq}
\expp{- \lambda t} \left( \e{\expp{\frac{\lambda}{k} h_L}} \right)^k + {n \choose m} \expp{- I(L)}
\hspace{3 cm} && \text{ Lemma } \ref{lem:mgfbd}
\\ & \leq
\expp{- \lambda t} \left( \expp{\frac{v(L, \frac{\lambda}{k})}{2} \frac{\lambda^2}{k^2}} \right)^k + {n \choose m} \expp{- I(L)}
\hspace{3 cm} && \text{ Lemma } \ref{lem:hmgfbd}
\\ & =
\expp{- \lambda t + \frac{v(L, \frac{\lambda}{k})}{2 k} \lambda^2} +  {n \choose m} \expp{- I(L)}.
\end{alignat*}
Choose $L = kt$.
To conclude the proof we only need to show that there are always choices for $\lambda$ which makes
$$ \expp{- \lambda t + \frac{v(kt, \frac{\lambda}{k})}{2 k} \lambda^2} \leq \expp{- \frac{k t^2}{2 v(kt, \beta \frac{I(kt)}{kt})}} + \expp{- \beta I(kt) \max(\frac{1}{2}, c(t, \beta, k))}$$.

We consider two cases:
\begin{enumerate}
    \item if $ \frac{t}{v \left( kt, \beta \frac{I(kt)}{kt} \right)} \leq \frac{\beta I(kt)}{kt}$
    choose $\lambda = \frac{k t}{v(kt, \beta \frac{I(kt)}{kt})}$,
    \hspace{.3 cm}
    so $\frac{\lambda}{k} = \frac{t}{v \left( kt, \beta \frac{I(kt)}{kt} \right)} \leq \frac{\beta I(kt)}{kt}$
    
    \item if $ \frac{t}{v \left( kt, \beta \frac{I(kt)}{kt} \right)} > \frac{\beta I(kt)}{kt}$
    choose $\lambda = \frac{\beta I(kt)}{t}$.
\end{enumerate}
Then, in case 1 since $\frac{\lambda}{k} \leq \frac{\beta I(L)}{L}$, we have $v(L, \frac{\lambda}{k}) \leq v(L , \frac{\beta I(kt)}{kt})$ (Note that $v(L, \cdot)$ is increasing in its second argument). Hence,

\begin{align*}
    - \lambda t + \frac{v(kt, \frac{\lambda}{k})}{2 k} \lambda^2 \leq  - \lambda t + \frac{v(kt, \frac{\beta I(kt)}{kt})}{2 k} \lambda^2 = - \frac{kt^2}{2 v(kt, \frac{\beta I(kt)}{kt})}.
\end{align*}

In the second case one just needs to substitute $\lambda$ to obtain
\begin{align*}
     - \lambda t + \frac{v(kt, \frac{\lambda}{k})}{2 k} \lambda^2 &= -{\beta I(kt)} + \frac{v(kt, \frac{\beta I(kt)}{kt} ) \beta^2 I(kt)^2 }{2 k t^2}
     \\ & =
     - \beta I(kt) \left(1 - \frac{v(kt, \frac{\beta I(kt)}{kt}) \beta I(kt)}{2kt^2} \right)
     \\ & =
      - \beta c(t, \beta, k) I(kt)
      \\ & = 
      - \beta \max(\frac{1}{2}, c(t, \beta, k)) I(kt).
\end{align*}

Note that since in this case $ \frac{t}{v \left( kt, \beta \frac{I(kt)}{kt} \right)} > \frac{\beta I(kt)}{kt}$, we have $c(t, \beta, k) > \frac{1}{2}$ so we have $\max(\frac{1}{2}, c(t, \beta, k)) = c(t, \beta, k)$.  The $\max$ operator controls this term when we are in the first case, so the upper bound does not blow up.

\end{proof}

\begin{remark}[Two regions of deviations] \label{rem: non-guassian boundary of upper bound}
    Inequality \eqref{eq:upper bound} reveals two different decay rates for the tail of $U_n$.  For small $t$s, the first term, i.e. $\expp{-\frac{k t^2}{2 v}}$, will be dominant, hence we observe Gaussian-like deviation.  This behavior has been studied already by CLT for U-statistics \cite{hoeffiding1948class}.  For larger $t$s, the last couple of terms on the right hand side of \eqref{eq:upper bound} will be dominant.  We call this region \textbf{large deviation} region.  Asymptotically, the sum of the last two terms decays like ${n \choose m} \expp{-I(kt)}$ for both subWeibull and polynomial tail kernels (see Section III of \cite{Bakhshizadeh2020} for detailed discussion).

    Inequality \eqref{eq:upper bound} denotes large deviation behavior whenever
    \begin{equation}
        \frac{k t^2}{v(kt, \beta )} \gg I(kt).
    \end{equation}
    For instance, when $I(kt) = \sqrt[\alpha]{kt}, \; \alpha \geq 1$ we have large deviation behavior for $ t \gg k^{- \frac{\alpha - 1}{2 \alpha - 1}} =  \lfloor \frac{n}{m} \rfloor^{- \frac{\alpha - 1}{2 \alpha - 1}}$.  This means the region of Gaussian deviation shrinks to $0$ as $n \to \infty$, when $\alpha > 1$.
\end{remark}

\subsubsection{Parameters of inequality \eqref{eq:upper bound}}
\label{sec: paramters of upper bound}

Theorem \ref{thm:general upper bound} bounds the tail of $U_n$ in terms of $k = \lfloor \frac{n}{m} \rfloor$ and the tail of kernel $h$.  The only terms of \eqref{eq:upper bound} that might seem unfamiliar are $c(t, \beta, k), v(kt, \beta \frac{I(k t)}{k t})$.  What typically happens in the asymptotic setting $n \to \infty$ is $c(t, \beta, k) \to 1, \; v(kt, \beta \frac{I(k t)}{k t}) \to Var(h)$.  Moreover, $\beta$ can be chosen arbitrarily close to $1$. Hence, for large $n$, one can think of upper bound \eqref{eq:upper bound} as $\expp{- \frac{k t^2}{2 Var(h)}} + (1 + {n \choose m}) \expp{- I(kt)}$.  For logarithmic $I(t)$ which corresponds to polynomial tail kernels, there are more restrictions on the constant $\beta$.  Nevertheless, for large deviation regime, the dominant term of \eqref{eq:upper bound} will still be ${n \choose m} \expp{-I(kt)}$.  This Section contains several statements that make the above claims precise. 

\begin{remark} \label{rem: v is almost variance}
Lemma 4 of \cite{Bakhshizadeh2020} states that $v(kt, \beta \frac{I(kt)}{kt}) \xrightarrow{k \to \infty} Var(h)$ in either of the following setups:
\begin{enumerate}
    \item $I(t) = c \sqrt[\alpha]{t}, \; \alpha > 1, \; \beta < 1$
    \item $I(t) = \gamma \log(t), \; \gamma > 2, \; \beta < 1 - \frac{2}{\gamma}$.
\end{enumerate}

Hence, for large values of $k$, one should be able to upper bound the first term of \eqref{eq:upper bound} with $\expp{-\frac{k t^2}{2 C Var(h)}}$, where $C < 1$, but can get arbitrarily close to $1$.

In the above, Case 1 includes all subWeibull variables, and Case 2 includes variables with polynomial tails and finite variances.

\end{remark}

\begin{remark} \label{rem: ct asymptotic}
    When $kt^2 \gg I(kt)$, and $v(kt, \beta \frac{I(kt)}{kt})$ is bounded, we have 
    \begin{equation} \label{eq: ct asymp}
        c(t, \beta, k) \xrightarrow{k \to \infty} 1.
    \end{equation}
    {This includes both cases of Remark \ref{rem: v is almost variance} with $t \gg k^{- \frac{\alpha - 1}{2 \alpha - 1}}$ and $t \gg \sqrt{\frac{\log k}{k}}$, respectively.}

    To verify \eqref{eq: ct asymp} it suffices to note that $c(t, \beta, k) = 1 - \frac{\beta}{2 t} \frac{I(kt)}{kt} v(kt, \beta \frac{I(kt)}{kt})$, $\; \frac{I(kt)}{kt} \xrightarrow{k \to \infty} 0$, and all other terms  in the definition of $c(t, \beta, k)$ remain bounded.
\end{remark}

While Remark \ref{rem: v is almost variance} provides asymptotic bounds for $v(kt, \beta \frac{I(kt)}{kt})$, one might need bounds for finite sample case to utilize Theorem \ref{thm:general upper bound}.  Next Lemma and Remark provide such bounds.

\begin{lemma} \label{lem: v is bounded for subweibull}
If $I(t) \geq c \sqrt[\alpha]{t}$ for some $\alpha \geq 1$, $\text{Var}(h) < \infty$, and $\beta < 1$ is fixed, then there is a fixed number $v < \infty$ such that for any $L > 1$ and  $ \eta \leq \beta \frac{I(L)}{L}$ we have $v(L, \eta) \leq v$.
\end{lemma}

\begin{proof}
    Since $\eta \leq \beta \frac{I(L)}{L}$, we have $v(L, \eta) \leq v(L, \beta \frac{I(L)}{L})$. Moreover, by Corollary 2 of \cite{Bakhshizadeh2020} we obtain

    \begin{align} \nonumber
        v(L, \beta \frac{I(L)}{L}) & \leq \e{h^2 \mathbf{1}(h \leq 0)} + \frac{\Gamma(2 \alpha + 1)}{((1 - \beta) c)^{2 \alpha}} + L^{\frac{1}{\alpha} - 1} \frac{\beta c \Gamma(3 \alpha + 1)}{3 ((1 - \beta) c)^{3 \alpha}}
        \\ & \leq \label{eq: v bound in subweibull case}
        \e{h^2} + \frac{\Gamma(2 \alpha + 1)}{((1 - \beta) c)^{2 \alpha}} + \frac{\beta c \Gamma(3 \alpha + 1)}{3 ((1 - \beta) c)^{3 \alpha}}.
    \end{align}
The right hand side of \eqref{eq: v bound in subweibull case} does not depend on $L$, hence it remains constant as $L \to \infty$.

Note that there is a slight change of variable for function $v(L, \eta)$ defined here and the one defined in \cite{Bakhshizadeh2020}, which of course has been taken into account in quoting Corollary 2 from \cite{Bakhshizadeh2020}. 

\end{proof}

\begin{remark} \label{rem: v upperbound for polynomial tail}
    If $I(t) \geq \gamma \log t$ with $\gamma > 2$, which includes kernels with polynomial tails and finite variances, Corollary 3 of \cite{Bakhshizadeh2020} yields
    \begin{equation}
        v(L, \beta \frac{I(L)}{L}) \leq C L^{2 - (1 - \beta) \gamma} \log L,
    \end{equation}
    for some $C$ independent of $L$.

    Hence, for $\beta < 1 - \frac{1}{\gamma}$, while $v(L, \beta \frac{I(L)}{L})$ can grow as $L \to \infty$, still the last two terms of \eqref{eq:upper bound} are the dominant terms of right hand side.  As discussed in \cite{Bakhshizadeh2020}, $\beta < 1 - \frac{1}{\gamma}$ is sufficient for obtaining sharp upper bounds for the deviation probability of the sum of iid variables with polynomial tails, i.e. U-statistics with order $m =1$. 
\end{remark}

\subsection{Lower bound} \label{sec: lowerbound}
Let $J(t)$ be the function defined in Definition \ref{def: I, J}, and $A_n = \intcc{- J^{-1}(\log 2n), J^{-1}(\log 2n)}^{n - 1}$, where $J^{-1}$ denotes generalized inverse of $J$ and $\intcc{\cdot, \cdot}$ is the closed interval of given limits.  
Define:
\begin{defn} \label{def: phi function}
    
\begin{equation} \label{def:phi}
    \varphi_n(X) \define \inf\limits_{(X_1, ..., X_{n - 1}) \in A_n} h(X_1, X_2, ..., X_{m- 1}, X).
\end{equation}
\end{defn}

Note that we force $X_1, ... X_{n-1}$ to be in $A_n$, but only use the first $m - 1$ in the argument of $\inf$.

Then we have:
\begin{lemma} \label{lem:lower bound based on phi}
Assume kernel $h$ has a finite variance. Then
\begin{align} \label{eq:lower bound}
    \p{U_n > t} \geq C \p{\varphi_n(X_n) \geq  \frac{nt}{m}},
\end{align}
where $C > 0$ is an absolute constant independent of $n$.
\end{lemma}

\begin{proof}
\begin{align*}
     \p{U_n > t} & \geq \p{U_{n-1}(X_1, ..., X_{n - 1}) \geq 0, \; \sum_{i_1< ...< i_{m -1}< n} h(X_{i_1}, X_{i_2}, ..., X_{i_{m-1}}, X_n) > {n \choose m} t}
     \\ & \geq
     \p{U_{n-1}(X_1, ..., X_{n - 1}) \geq 0, \; (X_1, ..., X_{n - 1}) \in A_n, {n - 1 \choose m - 1} \varphi_n(X_n) > {n \choose m} t }
     \\ & \geq
     \p{U_{n - 1}(X_1, ..., X_{n - 1}) \geq 0, \; \abs{X_i} \leq J^{-1}(\log 2n) \;\; \forall i \leq n -1} \p{\varphi_n(X_n) > \frac{n t}{m}}.
\end{align*}
Note that 
\begin{align*}
    \p{\abs{X_i} \leq J^{-1}(\log 2n) \;\; \forall i \leq n -1} &= \left( 1 - \p{\abs{X_i} > J^{-1}(\log 2n)} \right)^{n - 1}
    \\ & \geq
    \left( 1 - \expp{- J(J^{-1}(\log 2n))} \right) ^{n -1}
    \\ & \geq
    (1 - \frac{1}{2n})^{n - 1} \xrightarrow{n \to \infty} \frac{1}{\sqrt{\rm e}}.
\end{align*}
Moreover, $\p{U_{n - 1} \geq 0} \xrightarrow{n \to \infty} \frac{1}{2}$ by CLT for U-statistics \citep{hoeffiding1948class}.  Hence, for large enough $n$ we obtain:

\begin{equation*}
     \p{U_{n - 1}(X_1, ..., X_{n - 1}) \geq 0, \; \abs{X_i} \leq J^{-1}(\log 2n) \;\; \forall i \leq n -1} \geq 0.9 \left( \frac{1}{2} + \frac{1}{\sqrt{\rm e}} - 1 \right) > 0.
\end{equation*}
Choosing $C < 0.9 \left( \frac{1}{\sqrt{e}} - \frac{1}{2} \right)$, and small enough to cover all the finite cases before the above asymptotic become true concludes the proof.

\end{proof}

\begin{remark}
$A_n = \intcc{- J^{-1}(\log 2n), J^{-1}(\log 2n)}^{n - 1}$ in Lemma \ref{lem:lower bound based on phi} can be replaced with any sequence of events $A_n \subset \mathbb{R}^{n - 1}$ for which $\liminf\limits_{n \to \infty} \p{A_n} > \frac{1}{2}$.  Also, one can work with $\p{U_{n - 1} \geq - \epsilon, \varphi_n(X_n) > \frac{nt}{m} + \frac{\epsilon}{{n \choose m}}}$ to relax this condition to $\liminf\limits_{n \to \infty} \p{A_n} > 0$.
\end{remark}

\subsection{Large Deviation Principle} \label{sec: LDP}

In this Section, we show the upper bound \eqref{eq:upper bound} is asymptotically tight in certain cases.  The idea is to show the rate function for the large deviation of a U-statistic is asymptotically equivalent to the right hand side of \eqref{eq:upper bound}.  A trivial requirement for such sharpness to hold is functions $I(t), J(t)$ defined in Definition \ref{def: I, J} be asymptotically tight.  This is formalized in the next assumption.

\begin{assumption} \label{assump: I, J are asymptotically tight}
Suppose
    \begin{align*}
        \lim_{t \to \infty} \frac{- \log \p{h(X_1, ..., X_m) > t}}{I(t)} = \lim_{t \to \infty} \frac{- \log \p{\abs{X_i} > t}}{J(t)} = 1.        
    \end{align*}
\end{assumption}

Hereafter, we focus on subWeibull distributions.  The tail of such distribution is bounded by some Weibull distribution, hence, all of their moments are finite.  Nonetheless, the exponential moment is not finite.  Assumption \ref{assump: I(t) > root alpha t} encodes the class of heavy-tailed subWeibull random variables.

\begin{assumption} \label{assump: I(t) > root alpha t}
Assume there is $\alpha > 1$ and $c > 0$ such that $I(t) \geq c \sqrt[\alpha]{t}, \; \forall t > 0$.    
\end{assumption}

Although LDP is derived under Assumptions \ref{assump: I(t) > root alpha t} and \ref{assump: J is sup log} here, we think one can obtain similar results for distributions with polynomial tails, i.e. logarithmic $I(t), J(t)$, and finite second moments following footsteps of this Section and \citet{Bakhshizadeh2020}.  However, for the sake of brevity we do not include polynomial tails in the following Sections. 

Moreover, we need a lower bound for the deviation amount $t$ to make sure it is in the large deviation regime.  The below assumption provides such a bound.
\begin{assumption} \label{assump: kt^2 >> I(kt)}
    Suppose $\frac{k t^2}{I(kt)} \to \infty$, i.e. $k t^2 \gg I(kt)$, as $n \to \infty$.
\end{assumption}

\begin{remark}
    For $I(t) = c \sqrt[\alpha]{t}, \; \alpha > 1$, Assumption \ref{assump: kt^2 >> I(kt)} holds whenever $t \gg n^{- \frac{\alpha - 1}{2 \alpha - 1}}$.  This includes constant $t$ as well as converging to $0$ sequence $t_n$ as long it decays slower than $(\frac{1}{n})^{\frac{\alpha - 1}{2 \alpha - 1}}$.
\end{remark}

\begin{lemma} \label{lem: ldp}
Suppose Assumptions \ref{assump: centered kenrnel}, \ref{assump: I, J are asymptotically tight}, \ref{assump: I(t) > root alpha t}, and \ref{assump: kt^2 >> I(kt)} hold. For $\varphi_n(\cdot)$ defined in \eqref{def:phi} and $I(\cdot)$ in Theorem \ref{thm:general upper bound} if one has
\begin{equation} \label{eq:ldp condition}
    \lim_{n \to \infty} \frac{-\log \p{\varphi_n(X) \geq \frac{n}{m} t}}{ I(k t)} = 1,
\end{equation}
where $k = \lfloor \frac{n}{m} \rfloor$, then
\begin{equation*}
    \lim_{n \to \infty} \frac{- \log \p{U_n > t}}{I(k t)} = 1.
\end{equation*}

In other words, $I(k t)$ is the large deviation rate function for $U_n$.
\end{lemma}
We postpone proof of the above Lemma to Section \ref{secpf:  Lemma ref{lem: ldp}}.

Condition \eqref{eq:ldp condition} essentially says large values of $h(X_1, ..., X_m)$ are determined by large value of only one coordinate.  It can be proved for many commonly used kernels applied to heavy-tailed variables $X_i$.  Assuming  tail of $\abs{X_i}$ is a shifted sub-additive function, a usual property of heavy tails, we can prove \eqref{eq:ldp condition} for several kernels.

\begin{assumption} \label{assump: sub-additive J}
Suppose a constant shift on $J(t)$, defined in Definition \ref{def: I, J}, makes it sub-additive on non-negative Real numbers, i.e.
\begin{equation} \label{eq: sub-additive J}
    J(t_1 + t_2) \leq J(t_1) + J(t_2) + b, \quad \forall t_1, t_2 \geq 0,    
\end{equation}
where $b \in \mathbb{R}$ is an absolute constant.
\end{assumption}

\begin{remark}
Assumption \ref{assump: sub-additive J} is somehow posing heavy-tailed distribution requirement for random variable $X$.  While it does not directly state $J$ is a sublinear function, which is equivalent to the formal definition of a heavy-tailed distribution, it controls $J$'s growth to be equal to or slower than linear functions.  Lemma \ref{lem: concave functions are sub-additive} and Remarks \ref{rem: sub-additive distributions}, \ref{rem: other heavy-tailed dist are subadditive} denote most well-known heavy-tailed distributions can have shifted sub-additive tail function $J$.
\end{remark}

\begin{lemma} \label{lem: concave functions are sub-additive}
    If $J: \mathbb{R} \to \mathbb{R}$ is a concave function on non-negative Real numbers, then $J$ satisfies \eqref{eq: sub-additive J} with $b = - J(0)$.
\end{lemma}

\begin{proof}
   By concavity of $J$ we have
   \begin{equation*}
      J(t_1) =  J \left( \frac{t_1}{t_1 + t_2} (t_1 + t_2) + \frac{t_2}{t_1 + t_2} 0 \right) \geq \frac{t_1}{t_1 + t_2}J(t_1 + t_2) + \frac{t_2}{t_1 + t_2} J(0).
   \end{equation*}
   Similarly, we obtain $J(t_2) \geq \frac{t_2}{t_1 + t_2} J(t_1 + t_2) + \frac{t_1}{t_1 + t_2} J(0)$.  Summing up the above two inequalities shows 
   $$J(t_1 + t_2) \leq J(t_1) + J(t_2) - J(0).$$
\end{proof}

\begin{remark} \label{rem: sub-additive distributions}
    For the below distributions, one can find a function $J(t)$ as defined in Definition \ref{def: I, J} which is both asymptotically tight and shifted sub-additive, i.e. satisfies Assumptions \ref{assump: I, J are asymptotically tight}, \ref{assump: sub-additive J}.
    \begin{enumerate}
        \item Exponential
        \item $\abs{\mathcal{N}(0, 1)}^{\alpha}, \quad \alpha \geq 2$
        \item Log-Normal
        \item Weibull distribution with shape parameter $s \leq 1$
        \item Log Logistic
        \item Pareto 

    \end{enumerate}
\end{remark}
We postpone proof of Remark \ref{rem: sub-additive distributions} to Section \ref{secpf: Remark {rem: sub-additive distributions}}.

\begin{remark} \label{rem: other heavy-tailed dist are subadditive}
    In general, with simple modifications, we expect the tail function $J(t)$ for heavy-tailed distributions satisfy Assumption \ref{assump: sub-additive J}. This includes distributions that are not named in Remark \ref{rem: sub-additive distributions}.  Note that $J(t) = - \log \p{\abs{X} > t}$ is an increasing function that grows to infinity, and by heavy-tailed assumption is supposed to be sub-linear.  Hence, it is expected that $J(t)$ becomes a concave function after some point $t \geq T$.  Using the same technique we utilized in the proof of case 3, Log-Normal distribution, of Remark \ref{rem: sub-additive distributions}, one can define a function $J_2(t)$ which is equal to $J(t)$ on $[T, \infty)$, and linearly extends to $[0, T]$ such that it is less than $J(t)$ and remains concave on the whole non-negative Real numbers.  At this point, Lemma \ref{lem: concave functions are sub-additive} shows $J_2(t)$ should satisfy \eqref{eq: sub-additive J}. 
\end{remark}

\begin{assumption} \label{assump: J is sup log}
Suppose $J(t) \gg \log t$ as $t \to \infty$, i.e. $\lim\limits_{t \to \infty} \frac{\log t}{J(t)} = 0$.
\end{assumption}

\begin{lemma} \label{lem: kernels with ldp}
Under Assumptions  \ref{assump: I, J are asymptotically tight}, \ref{assump: sub-additive J}, \ref{assump: J is sup log} condition \eqref{eq:ldp condition} holds with any bounded $t$ in the following cases:
\begin{enumerate}
    \item $h(X, Y) = \abs{X - Y} - \e{\abs{X - Y}}$
    
    \item $h(X, Y) = (X - Y) ^2- \e{(X - Y)^2}$
    
    \item $h(X_1, ..., X_m) = \max (\abs{X_1}, ... ,\abs{X_m}) - \e{\max (\abs{X_1}, ... ,\abs{X_m}) }$ 
    
    \item $h(X, Y)=\frac{1}{2}(X^2 + Y^2) -\max{(X, Y)} - \e{X^2}  + \e{\max{(X, Y)}}$
\end{enumerate}

Hence, under extra Assumptions \ref{assump: I(t) > root alpha t}, \ref{assump: kt^2 >> I(kt)} Lemma \ref{lem: ldp} yields
\begin{equation*}
    \lim_{n \to \infty} \frac{ - \log \p{U_n > t}}{I(k t)} = 1,
\end{equation*}
for U-statistics constructed with the above kernels.
\end{lemma}
Proof of the above Lemma is postponed to Section \ref{secpf: lem: kernels with ldp}.

\begin{remark}
    The last kernel in Lemma \ref{lem: kernels with ldp} is related to $\omega^2$-statistics for the goodness of fit \cite{peaucelle2004efficiency}.
\end{remark}

\begin{remark} \label{rem: max kernel}
    Similar to case 3 of Lemma \ref{lem: kernels with ldp}, if one takes $J(t)$ to be the tail function of $X$ instead of $\abs{X}$, i.e. $\p{X > t} \lesssim \expp{-J(t)}$,  she can show condition \eqref{eq:ldp condition} also holds for the below kernel
    $$h(X_1, ..., X_m) = \max (X_1, ..., X_m) - \e{\max (X_1, ..., X_m)}.$$ 
\end{remark}

\subsection{Discussion on the necessity of condition \eqref{eq:ldp condition}}

Lemma \ref{lem: ldp} denotes the upper bound given in Theorem \ref{thm:general upper bound} is asymptotically sharp if \eqref{eq:ldp condition} holds.  Lemma \ref{lem: kernels with ldp} lists some common kernels of U-statistics for which \eqref{eq:ldp condition} holds.  One might ask if \eqref{eq:ldp condition} is necessary to obtain asymptotic sharpness and LDP.  Below, we study an example for which the bound given by Theorem \ref{thm:general upper bound} is not sharp.  This shows kernels need to satisfy certain conditions to have the same asymptotic as the right hand side of \eqref{eq:upper bound}.  Determining the necessary conditions for such kernels is not given in this work, but is an interesting question to be addressed in future studies.

Consider $m=2$ and $h(x,y)=xy$. Also, assume $X$ has a symmetric distribution around origin, and $J(t) = c t^{\alpha}, \; \alpha < 1$ (e.g $X \sim \text{Weibull}$ or $X = \N(0, 1)^{\frac{2}{\alpha}}$).  
In this case,
\begin{align*} 
     \p{XY > u} &\leq \p{\abs{X} > \sqrt{u} \text{ or } \abs{Y} > \sqrt{u}} \simeq 2\expp{-J(\sqrt{u})} = 2\expp{-c u^{\frac{\alpha}{2}}},
     \\ 
     \p{XY > u} & \geq  \p{X > \sqrt{u}} \p{Y > \sqrt{u}} \simeq \expp{-2J(\sqrt{u})} = \expp{-2 c u^{\frac{\alpha}{2}}},
\end{align*}
for large enough $u$.
Hence, for the tail function $I(u)$ we will have
\begin{equation} \label{eq: I growth for product kernel}
    C_1 u^{\frac{\alpha}{2}} \leq I(u) \leq C_2 u^{\frac{\alpha}{2}}, \quad C_1, C_2 > 0.
\end{equation}

If one directly applies Theorem \ref{thm:general upper bound}, she obtains

\begin{equation} \label{eq: xy kernel loose upper bound}
    \p{U_n > t} \leq \expp{-\frac{kt^2}{2 v}} + \expp{- \beta I(kt) \max(\frac{1}{2}, c(t, \beta, k))} + {n \choose 2} \expp{-I(kt)},
\end{equation}
where constant $v$ is given by Lemma \ref{lem: v is bounded for subweibull}, $\beta < 1$ is arbitrary for large $n$, and $c(t, \beta, k) \xrightarrow{n \to \infty} 1$.  The right hand side of \eqref{eq: xy kernel loose upper bound} is at least $\expp{-I(kt)} > \expp{-\frac{C_2}{2^{\frac{\alpha}{2}}} {n}^{\frac{\alpha}{2}}} = \expp{-C_2' n^{\frac{\alpha}{2}}}$.  This bound is loose since we will show $\p{U_n > t}$ decays like $\expp{-C_3 n^{\alpha}}$, and $n^{\alpha} \gg n^{\frac{\alpha}{2}}$ as $n \to \infty$.

Observe that
$$U_n=\frac{1}{n(n-1)}\left(\sum_{i=1}^n X_i\right)^2-\frac{1}{n(n-1)}\sum_{i=1}^n X_i^2 = \frac{n}{n - 1} S_n^2 - \frac{1}{n (n - 1)} T_n,$$
where $S_n = \frac{1}{n} \sum X_i, \; T_n = \sum X_i^2$.

Note that $S_n$ is an order $1$ U-statistic with kernel $h_2(x) = x$, and $U_n \leq \frac{n}{n - 1} S_n^2 \leq 2 S_n^2$. Also, for kernel $h_2$ the tail function $I_2(t) = J(t) - \log 2$.  Hence, by Theorem \ref{thm:general upper bound} we obtain

\begin{align} \nonumber
    \p{U_n > t} &\leq \p{S_n > \sqrt{\frac{t}{2}}} 
    \\ & \label{eq: xy kernel tight bound}
    \leq \expp{- \frac{n t}{4 v}} + C_3 \expp{-\beta \max(\frac{1}{2}, c(t, \beta, k)) J(n\sqrt{t/2})} + C_3 n \expp{-J(n \sqrt{t/2})}.
\end{align}

As discussed in Section \ref{sec: LDP}, the right hand side of \eqref{eq: xy kernel tight bound} decays like $\expp{ - J(n \sqrt{t/2})} = \expp{-C_3 n^{\alpha} t^{\frac{\alpha}{2}}}$, which is much smaller than \eqref{eq: xy kernel loose upper bound} when $n$ is large.

Indeed, we can show \eqref{eq: xy kernel tight bound} has the right order of decay.
Considering the event $E = (X_{n - 1}, X_{n} > \sqrt{n(n - 1) t} \text{ and } U_{n - 2} \geq 0)$ for which
 $$\p{U_n > t} \geq \p{E} \simeq \expp{-C_4 n^{\alpha} t^{\frac{\alpha}{2}}}.$$
   Therefore, one can show $n^{\alpha}$ is the correct speed for the large deviation decay of $U_n$.  In other words, there are constants $C_4, C_5 > 0$ such that for large enough $n$
 \begin{equation*}
     \expp{-C_4 (n \sqrt{t})^{\alpha}} \leq \p{U_n > t} \leq \expp{-C_5 (n \sqrt{t})^{\alpha}}.
 \end{equation*}

\begin{remark} \label{rem: thm can be sharp even for xy kernel}
    Note that $h(x, y) = xy$ is a degenerate unbounded kernel.  Both \citet{nikitin2001rough} and \citet{chakrabortty2018tail} claim sharp exponential bounds or large deviation limits for such kernels have not been addressed in works preceding them. As discussed in the current Section, while product kernel does not satisfy \eqref{eq:ldp condition}, a slight modification in the usage of Theorem \ref{thm:general upper bound} can still yield an exponential bound which is sharp up to a constant.  This shows the strength of Theorem \ref{thm:general upper bound}  even beyond scenarios in which sharpness has been shown through Lemma \ref{lem: ldp}.
\end{remark}

\section{Future works}

While we documented some important information about the concentration of U-statistics with heavy-tailed samples, there are questions that remained unanswered.  It seems addressing some of them takes only extra effort along the similar path of reasoning we used in this manuscript. We exclude such questions just for the sake of brevity which we believe helps to convey the main message of the current work better.  Other questions sound more challenging and may require different techniques to be addressed.  In this Section, we point out both types of questions and leave them for future studies.

Note that Lemmas \ref{lem: ldp} and \ref{lem: kernels with ldp} denote the upper bound \eqref{eq:upper bound} is sharp for certain U-statistics when $t$ is larger than the region of Gaussian decay. However, the first term of \eqref{eq:upper bound}, the Gaussian term, does not have a sharp constant in general.  Below Remark declares this fact.

\begin{remark}
Let $h_1(X_1) = \e{h(X_1, ..., X_m) \sVert X_1}$.  The asymptotic variance of $U_n$ is $\frac{m^2}{n} Var(h_1)$  \cite{hoeffiding1948class}, $v(kt, \beta \frac{I(kt)}{kt}) \to \text{Var}(h)$, and $Var(h) \geq m Var(h_1)$. In fact, $Var(h) - m Var(h_1) = Var(h - \sum\limits_{i \leq m} h_1(X_i))$.  This means 
$$\lim\limits_{n \to \infty} \frac{\frac{k t^2}{2 v(kt, \beta \frac{I(kt)}{kt})}}{ \frac{m t^2}{2 n \text{ Var }(U_n)}} < 1.$$

It would be interesting if one can improve \eqref{eq:upper bound} to have sharp constants on both Gaussian and heavy-tailed regions of deviation.  A direction that sounds promising to do so is to use Hoeffding decomposition \cite{hoeffiding1948class} and apply Theorem \ref{thm:general upper bound} for projections of $U_n$ individually.  

Another possible improvement is to extend results of Section \ref{sec: LDP} beyond kernels with subWeibull tails, i.e. when Assumption \ref{assump: I(t) > root alpha t} does not hold.  This already has been done for sums of iid samples, i.e. U-statistics of order $m  = 1$, in \citet{Bakhshizadeh2020}.  Moreover, Theorem \ref{thm:general upper bound} offers non-trivial bounds as long as $\text{Var}(X_i) < \infty$.  Assumption \ref{assump: I(t) > root alpha t} is used to remove all logarithmic terms when $n \to \infty$.  Taking such terms into account needs more effort, but does not change the spirit of our reasoning. 
Let $I(t) = \gamma \log t$ have logarithmic growth.  In the light of \eqref{eq:upper bound} $f(t) = \frac{I(kt) - \log {n \choose m}}{\log n} \simeq \frac{I(k t) - m \log n}{\log n}$ seems a reasonable rate function for LDP of $U_n$ with speed $b_n = \log n$.

Moreover, condition \eqref{eq:ldp condition} is only a sufficient condition for the sharpness of inequality \eqref{eq:upper bound}.  It is interesting to determine necessary conditions for sharpness of Theorem \ref{thm:general upper bound} in the sense of rough logarithmic limits.  In addition, developing sharp upper bounds when such conditions do not hold can be a good direction to extend results of the current paper.

Finally, perhaps the most important message of this paper is one can extend concentration results developed for subGaussian variables to distribution with heavier tails simply by truncation technique and precise tuning of the truncation level.  While this technique is applied for U-statistics here, the question is to what other self-averaging processes we can apply the same and obtain sharp concentration.  We hope this work motivates future studies to obtain upper bounds with sharp constants to a larger class self-averaging processes.

\end{remark}

\section*{Acknowledgment}
The author is thankful for inspiring discussions he has had with Dr. Nabarun Deb during the development of this work. 

\section{Proofs}
This Section includes the proofs of statements in the previous Sections.  First, let recall the following Lemma from \citet{Bakhshizadeh2020} which turns out to be useful in the calculation of logarithmic limits.

\begin{lemma}[Lemma 5 of \cite{Bakhshizadeh2020}] \label{lem: log of sum limits}
    Let $a_n, b_n$ and $c_n$ be sequences of positive numbers such that
    \begin{equation*}
        \lim_{n \to \infty} \frac{\log a_n}{c_n} = a, \;
        \lim_{n \to \infty} \frac{\log b_n}{c_n} = b, 
        \;
        \lim_{n \to \infty} c_n = \infty.
    \end{equation*}
    Then
    \begin{equation*}
        \lim_{n \to \infty} \frac{\log (a_n + b_n)}{c_n} = \max \cbr{a, b}.
    \end{equation*}
\end{lemma}

\subsection{Proof of Lemma \ref{lem: ldp}} \label{secpf:  Lemma ref{lem: ldp}}

\begin{proof}

By Lemma \ref{lem: v is bounded for subweibull}, for any $\beta < 1$, we have $v(kt, \beta \frac{I(k t)}{kt}) < v$, where $v$ is a constant independent of $k$.  Therefore, we obtain the following
\begin{equation*}
     \lim_{k \to \infty} \frac{\frac{- k t^2}{2 v(kt, \beta \frac{I(kt)}{kt})}}{I(kt)} \leq \lim_{k \to \infty}\frac{\frac{-k t^2}{2 v}}{I(kt)} = -\infty,
\end{equation*}
because $\frac{kt}{I(kt)} \to \infty$ as $kt \to \infty$.
Moreover, by Remark \ref{rem: ct asymptotic}, $c(t, \beta, k) \to 1$, hence,
\begin{equation*}
    \lim_{k \to \infty} \frac{- \beta I(kt) \max(\frac{1}{2}, c(t, \beta, k))}{I(kt)} = -\beta, \quad \lim_{k \to \infty} \frac{ \log {n \choose m} - I(kt)}{I(kt)} = -1.
\end{equation*}

Having inequality \eqref{eq:upper bound} and the above equations, Lemma \ref{lem: log of sum limits} yields $\lim\limits_{n \to \infty} \frac{ \log \p{U_n > t}}{I(k t)} \leq - \beta, \; \forall \beta < 1$.  Multiplying by $-1$ and taking supremum over $\beta < 1$ implies
\begin{equation}
    \lim_{k \to \infty} \frac{- \log \p{U_n > t}}{I(k t)} \geq 1.
\end{equation}

On the other hand, Lemma \ref{lem:lower bound based on phi} yields
$
    \lim\limits_{n \to \infty} \frac{- \log \p{U_n > t}}{-\log \p{\varphi_n(X) > \frac{n t}{m}}} \leq 1,
$
so by \eqref{eq:ldp condition} we obtain
\begin{equation*}
    \lim_{n \to \infty} \frac{- \log \p{U_n > t}}{I(kt)} \leq 1,
\end{equation*}

\end{proof}

\subsection{Proof Remark \ref{rem: sub-additive distributions}} \label{secpf: Remark {rem: sub-additive distributions}}

\begin{proof}
    \hfill
    \begin{enumerate}     
        \item 
        If $X \sim \text{Exp}(\lambda)$ we have $\p{X > t} = \expp{-\lambda t}$, so $J(t) = -\log \p{X > t} = \lambda t$ is a linear function which satisfies \eqref{eq: sub-additive J} with equality and $b = 0$.

        \item Let $X = \abs{Z}^\alpha, \; Z \sim \N(0, 1)$.  Observe that
        \begin{equation*}
            \p{X > t} = \p{\abs{Z} > \sqrt[\alpha]{t}} \leq 2\expp{-\frac{1}{2} t^{\frac{2}{\alpha}}}
        \end{equation*}
        Hence, setting $J(t) = \frac{\sqrt[\frac{\alpha}{2}]{t}}{2} - \log 2 $ gives us a tail upperbound as in \eqref{eq: J def}.  Utilizing LDP rate function of the Normal distribution shows $J(t)$ is asymptotically tight too.  Also, for $\alpha \geq 2, \; J(t)$ is a concave function, hence by Lemma \ref{lem: concave functions are sub-additive} satisfies Assumption \ref{assump: sub-additive J}. 

        \item
        Let $X = \expp{Z}, \; Z \sim \N(0, 1)$.  Note that
        \begin{equation*}
            \p{X > t} = \p{Z > \log t} \leq \expp{-J(t)},
        \end{equation*}
        where
        $    J(t) = \begin{cases}
                0 &  0 \leq t \leq 1
                \\
                {\frac{1}{2} \log^2 t} & t > 1    
            \end{cases}
            $.  Similar to the previous case, asymptotic tightness of $J$ comes from rate function of the Normal distribution.  Instead of showing \eqref{eq: sub-additive J} for $J(t)$, we replace it with a concave asymptotically tight lower bound $J_2(t)$.  This useful technique can be applied to several other cases as well.  Let
            \begin{equation*}
                J_2(t) = \begin{cases}
                    \frac{1}{\rm e}(t - {\rm e}) + \frac{1}{2} & 0 \leq t \leq {\rm e}
                    \\
                    \frac{1}{2} \log^2 t & t > {\rm e}
                \end{cases}.
            \end{equation*}
            Then, $J(t) \geq J_2(t) \; \forall t \geq 0 $, i.e. $\p{X > t} \leq \expp{-J_2(t)}$, and $J_2$ is a concave function on $\R^{\geq 0}$, hence by Lemma \ref{lem: concave functions are sub-additive} satisfies \eqref{eq: sub-additive J}.

            To verify above claims, one only needs to note that $J(t)$ is a convex function on $[0, {\rm e}]$ and is a concave function on $[{\rm e}, \infty)$.  $J_2(t)$ on $[0, \rm e]$ is simply the linear approximation of $J(t)$ at $t = {\rm e}$ to make it concave everywhere.

            \begin{figure}[h!]
                \centering
                \includegraphics{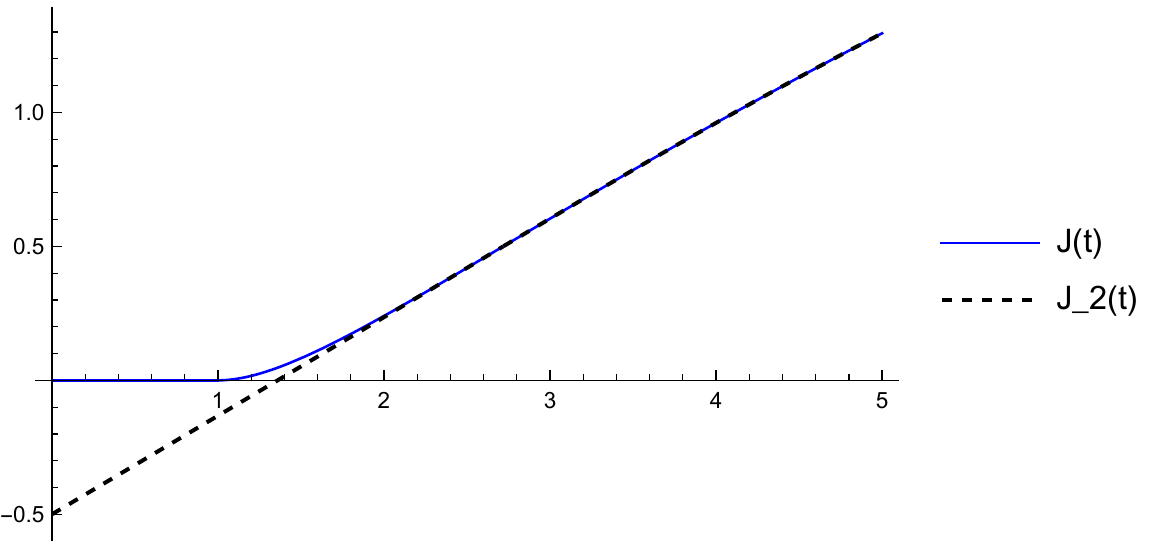}
                \caption{$J(t)$ and $J_2(t)$}
                \label{fig: j vs j2}
            \end{figure}
        
        \item If $X \sim \text{Weibull}(\lambda, s)$, then $\p{\abs{X} > t} = \expp{- (\frac{t}{\lambda})^s}$.  Hence, one can take the trivial tail bound $J(t) = - \log \p{X > t} = \left( \frac{t}{\lambda} \right)^s$. This function satisfies $J(0) = 0$, and is concave for $s \leq 1$.  Hence, Lemma \ref{lem: concave functions are sub-additive} yeilds $J(t_1 + t_2) \leq J(t_1) + J(t_2)$.

        \item Let $X \sim \text{Log-logistic}(\alpha, \beta), \; \alpha, \beta > 0$
        $$J(t) = - \log \p{\abs{X} > t} = - \log(1 - \frac{1}{1 + (x / \alpha)^{- \beta}}) = \log (1 + (x / \alpha)^{\beta}).$$

        Note that 
        $$ 2^{\beta} (1 + x^{\beta}) (1 + y^{\beta}) \geq 2^{\beta} + (2 \max\{x, y\})^{\beta} \geq 1 +  (x + y)^{\beta}, \quad \forall x, y \geq 1.$$
        therefore
        \begin{equation*}
            2^{\beta} (1 + (t_1 / \alpha)^{\beta})(1 + (t_2 / \alpha)^{\beta}) \geq (1 + ((t_1 + t_2)/ \alpha)^{\beta}).
        \end{equation*}
        Applying $\log$ to the above inequality yields
        $J(t_1) + J(t_2) + \beta \log 2 \geq J(t_1 + t_2)$.

        \item Let $X \sim \text{Pareto}(x_m, \alpha)$. Then,
        $$J(t) = - \log \p{X > t} = \alpha \log \frac{t}{x_m}, \quad t > x_m,$$ is a concave function on the support of $X$, i.e. $[x_m, \infty)$.  Similar to the Log-Normal case above, one can linearly extend $J(t)$ to a concave function on $[0, \infty)$ and utilize Lemma \ref{lem: concave functions are sub-additive} to verify $J(t)$ is shifted sub-additive.
    \end{enumerate}
\end{proof}

\subsection{Proof of Lemma \ref{lem: kernels with ldp}} \label{secpf: lem: kernels with ldp}

\begin{proof}

\emph{1.} The strategy to show \eqref{eq:ldp condition} is to show
$ \lim \frac{- \log \p{\varphi_n(X) > \frac{n}{m} t}}{J(kt)} = \lim \frac{I(kt)}{J(kt)} = 1$ as $n \to \infty$.
Let us write $c \define \e{\abs{X-Y}}$. Note that 

\begin{align*}
\varphi_n(X) &= \inf_{\abs{y} \leq J^{-1}(\log{2n})} \abs{X - y} -c 
\\ &
=\begin{cases}
-c & \mbox{if}\ \abs{X} \leq J^{-1}(\log{2n})
\\ \abs{X} -J^{-1}(\log{2n})-c & \mbox{if}\ \abs{X} >J^{-1}(\log{2n}).
\end{cases}.
\end{align*}
With the above display in mind, observe that
\begin{align*}
    \bP(\varphi_n(X)> \frac{n}{m} t)&=\p{\abs{X} >  \frac{n}{m} t + J^{-1}(\log 2n) + c}.
\end{align*}
 Note that, by Assumption \ref{assump: I, J are asymptotically tight} we have:
$$\lim\limits_{n\to\infty}\frac{-\log\bP(\abs{X}> \frac{n}{m} t + J^{-1}(\log{2n})+c)}{J(\frac{n}{m} t + J^{-1}(\log 2n) + c)}=1.$$

Also, note that for large $n$,  $J(k t) \leq J(\frac{n}{m} t + J^{-1}(\log 2n) + c) \leq J(kt + t + J^{-1}(\log 4 k m ) + c)$.
Using  \eqref{eq: lem J asymp 1} from Lemma \ref{lem: J asymptotics does not change with log} with $u = 4 k m, \; c_1 = \frac{t}{4m}, \; c_2 = c + t$ we obtain 
$  \lim\limits_{n \to \infty} \frac{J(k t + t + J^{-1}(\log 4 k m) + c)}{J(k t)} =  1 $,
therefore 

\begin{align}\label{eq:lub1}
        \lim\limits_{n\to\infty}\frac{-\log{\bP(\varphi_n(X)> \frac{n}{m} t)}}{J(kt)} =  \lim\limits_{n\to\infty}\frac{J(\frac{n}{m} t + J^{-1}(\log 2n) + c)}{J(kt)} = 1.
\end{align}
Next, we try to approximate the term $I(kt)$ from~\eqref{eq:ldp condition}.
Towards this direction, we begin by observing that
\begin{equation*}
\bP(|X_1-X_2|\geq c+kt) \geq  \bP(X_2 \leq  c')\bP(X_1\geq c' + c + kt),
\end{equation*}
for some $c' \in \mathbb{R}$ such that $\p{X_2 \leq c'} > 0$.
Consequently,
\begin{equation*}
\limsup\limits_{n\to\infty}\frac{-\log{\bP(h(X_1,X_2)\geq  kt)}}{J(kt + c' + c)}\leq 1,
\end{equation*}
Hence,
\begin{equation}\label{eq:lub1.2}
\limsup\limits_{n \to\infty} \frac{I(kt)}{J(kt + c' + c)} = \limsup\limits_{n\to\infty} \frac{I(kt)}{J(kt)} \leq 1.
\end{equation}
We used Lemma \ref{lem: J asymptotics does not change with log} to drop $c' + c$.

For the other direction, we need to establish an upper bound for $\p{h(X_1, X_2) > kt}$.  Let $u > 0$, then

\begin{align*}
    \bP(|X_1-X_2|\geq u) &\leq \p{\abs{X_1} > u} + \p{\abs{X_2} > u} + \p{\abs{X_1 - X_2} \geq u, \; \abs{X_1}, \abs{X_2} \leq u}
    \\ & \leq
    2 \expp{- J(u)} + 2 \p{X_1 > X_2 + u, \; \abs{X_1}, \abs{X_2} \leq u}
    \\ & \leq
    2 \expp{- J(u)} + 2 \sum_{i = 0}^{\lceil u \rceil} \p{- \frac{i}{\lceil u \rceil} u \leq X_2 \leq - \frac{i - 1}{\lceil u \rceil} u} \p{X_1 \geq (1 - \frac{i}{{\lceil u \rceil}}) u}
    \\ & \leq
    2 \expp{- J(u)} + 2 \sum_{i = 0}^{\lceil u \rceil} \expp{- J \paranth{ \frac{i - 1}{\lceil u \rceil} u} - J \paranth{ ( 1 - \frac{i}{{\lceil u \rceil}}) u}}
    \\ & \overset{*}{\leq}
    2 \expp{- J(u)} + 2 \sum_{i = 0}^{\lceil u \rceil} \expp{- J \paranth{u - \frac{u}{\lceil u \rceil}} + b}
    \\ & \leq
    2 \expp{- J(u)} +  2 (u + 2) \expp{b} \expp{- J(u - 1)}
    \\ & \leq
    3 {\rm e}^{b} u \expp{-J(u - 1)}, \quad \text{for } u \geq 2 {\rm e}^{-b} + 4.
\end{align*}

To obtain inequality marked by $*$, we used Assumption \ref{assump: sub-additive J}.
Taking $- \log$ of above inequality we get
\begin{equation} \label{eq: liminf -log abs}
    \liminf_{u \to \infty} \frac{- \log \p{\abs{X_1 - X_2} > u}}{- \log 3u - b + J(u - 1)} \geq 1.
\end{equation}
Set $u = kt + c$.  Since, by Lemma \ref{lem: J asymptotics does not change with log}, $\lim\limits_{k \to \infty} \frac{- \log 3(kt + c) - b + J(kt + c - 1)}{J(kt)} = 1$,  we obtain:
\begin{equation} \label{eq:abs I,J upbd}
    \liminf_{u \to \infty} \frac{I(kt)}{J(kt)} \geq 1.
\end{equation}

Equations \eqref{eq:lub1}, \eqref{eq:lub1.2}, and \eqref{eq:abs I,J upbd} yield \eqref{eq:ldp condition}.

\emph{2.} Once again, we set $c \define \E(X-Y)^2$ and note that, in this case,
\begin{align*}
\varphi_n(X)&=\inf_{\abs{y} \leq J^{-1}(\log{2n})} (X - y)^2 - c \\ &=
\begin{cases} 
-c & \mbox{if}\ \abs{X} \leq J^{-1}(\log{2n}), 
\\
(\abs{X}-J^{-1}(\log{2n}))^2-c &\mbox{if}\ \abs{X} > J^{-1}(\log{2n}).
\end{cases}
\end{align*}
With the above display in mind, observe that
\begin{align*}
    \bP(\varphi_n(X)> \frac{n}{m} t)&=\bP(\abs{X} >J^{-1}(\log{2n}) + \sqrt{c+ \frac{n}{m} t}).   
\end{align*}

Note that for large enough $n$

$$ \sqrt{kt} \leq  J^{-1}(\log{2n}) + \sqrt{c + \frac{n}{m} t}  \leq 
J^{-1}(\log{2n}) + \sqrt{\abs{c}} + \sqrt{t} + \sqrt{kt}.$$

Hence, by Assumptions \ref{assump: sub-additive J}, \ref{assump: J is sup log} we obtain
$$\lim\limits_{n\to\infty}\frac{J( J^{-1}(\log{2n}) + \sqrt{c + \frac{n}{m} t} )}{J(\sqrt{kt})}=1.$$
(see proof of Lemma \ref{lem: J asymptotics does not change with log} for details.)

We therefore have
\begin{align}\label{eq:lub4}
    \lim\limits_{n\to\infty}\frac{-\log{\bP(\varphi_n(X)> \frac{n}{m} t)}}{J(\sqrt{kt})}=1.
\end{align}
Next, we try to approximate the term $I(kt)$ from~\eqref{eq:ldp condition}.
Note that
\begin{equation}\label{eq:lub5}
\bP((X_1-X_2)^2\geq c+kt) = \bP(\abs{X_1 - X_2} \geq \sqrt{c+kt})\geq \bP(\abs{X_2} \leq c')\bP(\abs{X_1} \geq c' + \sqrt{c+kt}),
\end{equation}
for a constant $c'$ such that $\p{\abs{X_2} \leq c'} > 0$.
Consequently,
\begin{equation*}
\limsup\limits_{n\to\infty}\frac{-\log{\bP((X_1-X_2)^2\geq c+kt)}}{J(c' + \sqrt{c + kt})}\leq 1,
\end{equation*}
which in turn yields
\begin{equation}\label{eq:lub6}
\limsup\limits_{n\to\infty} \frac{I(kt)}{J(\sqrt{kt})}\leq 1.
\end{equation}
We utilized Lemma \ref{lem: J asymptotics does not change with log} to drop $c, c'$ from denominator in limits.

For the other direction, we need to establish an upper bound for the left hand side of \eqref{eq:lub5}. 
As proved in \eqref{eq: liminf -log abs}, with $u = \sqrt{c + kt}$, we have
\begin{equation} \label{eq: liminf abs sqrt}
    \liminf_{k \to \infty} \frac{- \log \p{\abs{X_1 - X_2} > \sqrt{c + kt}}}{- \log 3 \sqrt{c + kt} - b + J(\sqrt{c + kt} - 1)} \geq 1.
\end{equation}
Since $\frac{\log\sqrt{c + kt}}{J(\sqrt{kt})} \to 0, \; \frac{J(\sqrt{c + kt} - 1)}{J (\sqrt{kt})} \to 1$ as $k \to \infty$, from \eqref{eq: liminf abs sqrt} we obtain
\begin{equation}
    \liminf_{k \to \infty} \frac{- \log \p{h(X_1, X_2) \geq kt}}{J(\sqrt{k t})} \geq 1.
\end{equation}
This completes the proof.

\vspace{0.1in}

\emph{3.} We write $c=\E \max (\abs{X_1}, ... ,\abs{X_m})$. Note that
\begin{align*}
\varphi_n(X)= \abs{X} - c.
\end{align*}

Hence, $\p{\varphi_n(X) \geq \frac{n}{m} t} = \p{\abs{X} \geq \frac{n}{m} t + c}$, and similar to the above cases we can show:

\begin{equation}
    \lim_{n \to \infty} \frac{- \log \p{\varphi_n(X) \geq \frac{n}{m} t}}{J(kt)} = \lim_{n \to \infty} \frac{I(k t)}{J(kt)} = \lim_{n \to \infty} \frac{- \log \p{\varphi_n(X) \geq \frac{n}{m} t}}{I(kt)} =  1.
\end{equation}

\emph{4.} Assume $n$ is large enough so $1 \leq J^{-1}(\log 2n)$.  Calling $\e{X^2} - \e{\max\cbr{X, Y}} = c$ we have
\begin{equation*}
    \varphi_n(X) = 
    \begin{cases} 
    \frac{1}{2} X^2 - X - c, &  \mbox{if} \ X > \frac{1}{2},
    \\
    \frac{1}{2} X^2 - \frac{1}{2} - c, & \mbox{if} \ X \leq \frac{1}{2}.
    \end{cases}
\end{equation*}
Therefore, 
\begin{align*}
    \p{\varphi_n(X) > \frac{n}{m} t} & = \p{\frac{X^2}{2} - X > \frac{n}{m} t + c, X > 0} + \p{X \leq - \sqrt{\frac{2n}{m} t + 2c + 1}}
    \\ & =
    \p{X > 1 + \sqrt{\frac{2 n}{m} t + 2 c + 1}} + \p{X \leq - \sqrt{ \frac{2 n}{m} t + 2c + 1}}.
\end{align*}
Hence,
\begin{equation*}
    \p{\abs{X} > \sqrt{\frac{2n}{m} t + 2 c + 1} + 1} \leq \p{\varphi_n(X) > \frac{n}{m} t} \leq \p{\abs{X} > \sqrt{\frac{2n}{m} t + 2 c + 1}}.
\end{equation*}
Similar to the previous cases above, by Lemma \ref{lem: J asymptotics does not change with log}, we can then show
\begin{equation}
    \lim_{n \to \infty} \frac{- \log \p{\varphi_n(X) > \frac{n}{m} t}}{J \left( \sqrt{\frac{2n}{m} t + 2c + 1} + 1 \right)} =  \lim_{n \to \infty} \frac{- \log \p{\varphi_n(X) > \frac{n}{m} t}}{J \left( \sqrt{\frac{2n}{m} t + 2c + 1} \right)} = \lim_{n \to \infty} \frac{- \log \p{\varphi_n(X) > \frac{n}{m} t}}{J(\sqrt{2 k t})} = 1.
\end{equation}

Moreover, since $h(X, Y) \geq \min (\frac{1}{2} X^2 - \frac{1}{2}, \frac{1}{2} X^2 - X) - c$, we obtain

\begin{align*}
    \p{h(X, Y) > kt } &\geq \p{\min (\frac{1}{2} X^2 - \frac{1}{2}, \frac{1}{2} X^2 - X) > kt +  c} 
    \\ &\geq
    \p{\abs{X} > \sqrt{2 k t + 2 c + 1} + 1}.
\end{align*}

Hence,
\begin{equation}
    \limsup_{n \to \infty} \frac{- \log \p{h(X, Y) > k t}}{J(\sqrt{2 k t + 2 c + 1} + 1)} = \limsup_{n \to \infty} \frac{I(k t)}{J(\sqrt{2 k t} )} \leq 1.
\end{equation}
Again, for dropping extra constants from the denominator we used Lemma \ref{lem: J asymptotics does not change with log}.

Furthermore,
$h(X, Y) \leq \frac{1}{2} X^2 + \abs{X} + \frac{1}{2} Y^2 + \abs{Y} - c$, and $\frac{1}{2} X^2 + \abs{X} > u \iff \abs{X} > - 1 + \sqrt{2 u + 1}, \; \forall u \geq 0$.  Hence,

\begin{align*}
    \p{h(X, Y) > k t} & \leq \p{\frac{1}{2} X^2 + \abs{X} + \frac{1}{2} Y^2 + \abs{Y} > k t + c}
    \\ & \leq
    \p{\frac{1}{2} X^2 + \abs{X} > k t  + c} + \p{\frac{1}{2} Y^2 + \abs{Y} > k t  + c}
    \\ & \quad
    + \sum_{i = 1}^{\lceil kt + c \rceil} \p{\frac{i - 1}{{\lceil kt + c \rceil}} (kt + c) \leq \frac{1}{2} X^2 + \abs{X} \leq \frac{i}{{\lceil kt + c \rceil}} (kt + c)} \times
    \\ & \hspace{2 cm}
    \p{\frac{1}{2} Y^2 + \abs{Y} > (1 - \frac{i}{{\lceil kt + c \rceil}}) (kt + c)}
    \\ & \leq
    2 \p{\abs{X} > -1 + \sqrt{2 kt + 2c + 1}} 
    \\ & \quad
    + \sum_{i = 1}^{\lceil kt + c \rceil} \p{\abs{X} \geq -1 + \sqrt{\frac{2(i - 1)}{{\lceil kt + c \rceil}} (kt + c) + 1}} \times
    \\ & \hspace{2 cm}
    \p {\abs{Y} > -1 + \sqrt{2(1 - \frac{i}{{\lceil kt + c \rceil}}) (kt + c) + 1}}
    \\ & \leq
    2 \expp{-J(-1 + \sqrt{2 k t + 2c + 1})} 
    \\ & \quad
    + \sum_{i = 1}^{\lceil kt + c \rceil} \expp{-J (-1 +  \sqrt{\frac{2(i - 1)}{{\lceil kt + c \rceil}} (kt + c) + 1}) - J( -1 + \sqrt{2(1 - \frac{i}{{\lceil kt + c \rceil}}) (kt + c) + 1})}
    \\ & \overset{*}{\leq}
    2 \expp{-J(-1 + \sqrt{2 k t + 2c + 1})}  + \sum_{i = 1}^{\lceil kt + c \rceil} {\rm e}^{b} \expp{-J(-2 + \sqrt{2(1 - \frac{1}{\lceil kt + c \rceil}) (k t + c) + 2} )}
    \\ & =
    2 \expp{-J(-1 + \sqrt{2 k t + 2c + 1})} + {\lceil kt + c \rceil} {\rm e}^b \expp{-J(-2 + \sqrt{2(1 - \frac{1}{\lceil kt + c \rceil}) (k t + c) + 2})}.
\end{align*}
To obtain inequality marked by $*$ we used Assumption \ref{assump: sub-additive J} and the fact that $\sqrt{x + y} \leq \sqrt{x} + \sqrt{y}$ for non-negative $x, y$.

Hence, by Lemma \ref{lem: log lim for w2 kernel} we obtain
\begin{equation}
    \liminf_{n \to \infty} \frac{I(kt)}{ J(\sqrt{2 k t})} \geq 1,
\end{equation}
which concludes the proof.
\end{proof}

\begin{lemma} \label{lem: J asymptotics does not change with log}
If tail function $J(t)$, defined in \eqref{eq: J def}, satisfy Assumptions \ref{assump: sub-additive J}, \ref{assump: J is sup log}, then
for any $c_1 > 0$ and  $c_2, c_3 \in \mathbb{R}$ we have

\begin{equation}  \label{eq: lem J asymp 0}
    \lim_{u \to \infty} \frac{J(u + c_2)}{J(u)} = 1
\end{equation}

\begin{equation} \label{eq: lem J asymp 1}
     \lim_{u \to \infty} \frac{J(c_1 u + J^{-1}(\log u) + c_2 )}{J(c_1 u)} = 1
\end{equation}

\begin{equation} \label{eq: lem J asymp 2}
    \lim_{u \to \infty} \frac{J(c_2 + \sqrt{u + c}_3)}{J(\sqrt{u})} = 1
\end{equation}
\end{lemma}

\begin{proof}
To prove \eqref{eq: lem J asymp 0}, one only need to note that by Assumption \ref{assump: sub-additive J} we have
\begin{equation}
     J(u) - J(\abs{c_2}) - b \leq J(u - \abs{c_2}) \leq J(u + c_2) \leq J(u + \abs{c_2}) \leq J(u) + J(\abs{c_2}) + b,
\end{equation}
and that $J$ is a non-decreasing function, and $J(u) \xrightarrow{u \to \infty} \infty$. 

To show \eqref{eq: lem J asymp 1} observe that $J^{-1}(\log u) + c_2 > 0$ for large enough $u$.  
By Assumption \ref{assump: sub-additive J} we obtain
\begin{align*}
    J(c_1 u) \leq J(c_1 u + J^{-1}(\log u) + c_2 ) & \leq J(c_1 u) + \log u + J(\abs{c_2}) + 2b.
\end{align*}
Note that $\lim \frac{\log u}{J(u)} = \lim \frac{C}{J(u)} = 0$ as $u \to \infty$, so dividing by $J(c_1 u)$ and taking limit of the above inequalities yields \eqref{eq: lem J asymp 1}.

For the third part observe that when $u$ is large enough we have
\begin{equation*}
    \sqrt{u} - \sqrt{\abs{c_3}} - \abs{c_2} \leq c_2 + \sqrt{u + c_3} \leq \sqrt{u} + \sqrt{\abs{c_3}} + \abs{c_2}.
\end{equation*}
Then, since $J$ is non-decreasing and satisfies Assumption \ref{assump: sub-additive J} we obtain
\begin{equation*}
        J(\sqrt{u}) - J(\sqrt{\abs{c_3}} + \abs{c_2}) - b \leq J(c_2 + \sqrt{u + c_3}) \leq J(\sqrt{u}) + J(\sqrt{\abs{c_3}} + \abs{c_2}) + b.
\end{equation*}
Once again, dividing by $J(\sqrt{u})$ and taking $u \to \infty$ yields \eqref{eq: lem J asymp 2}.
\end{proof}

\begin{lemma} \label{lem: log lim for w2 kernel}
Under Assumptions \ref{assump: sub-additive J}, \ref{assump: J is sup log} we have
    \begin{equation*}
        \lim_{k \to \infty} \frac{ - \log \left( 2 \expp{-J(-1 + \sqrt{2 k t + 2c + 1})} + {\lceil kt + c \rceil} {\rm e}^b \expp{-J(-2 + \sqrt{2(1 - \frac{1}{\lceil kt + c \rceil}) (k t + c) + 2})} \right)} {J(\sqrt{2 k t})} = 1
    \end{equation*}
\end{lemma}

\begin{proof}
   Given Lemma \ref{lem: log of sum limits}, it suffices to show that
   \begin{align}
        & \lim_{k \to \infty} \frac{ - \log \left( 2 \expp{-J(-1 + \sqrt{2 k t + 2c + 1})} \right)}{ J(\sqrt{2 k t})} = 1
        \\ &
         \lim_{k \to \infty} \frac{- b - \log \left( {\lceil kt + c \rceil} \expp{-J(-2 + \sqrt{2(1 - \frac{1}{\lceil kt + c \rceil}) (k t + c) + 2})} \right)} {J(\sqrt{2 k t})} 
        = 1.
   \end{align}

   Note that
   \begin{align*}
       \lim_{k \to \infty} \frac{ - \log \left( 2 \expp{-J(-1 + \sqrt{2 k t + 2c + 1})} \right)}{ J(\sqrt{2 k t})}
       & =
       \lim_{k \to \infty} \frac{J(-1 + \sqrt{2 k t + 2c + 1})}{J(\sqrt{2 k t})} = 1,
   \end{align*}
   by Lemma \ref{lem: J asymptotics does not change with log}.  

   Moreover,
   \begin{align*}
       &\lim_{k \to \infty} \frac{-b - \log \left( {\lceil kt + c \rceil} \expp{-J(-2 + \sqrt{2(1 - \frac{1}{\lceil kt + c \rceil}) (k t + c) + 2})} \right)} {J(\sqrt{2 k t})} 
       \\ & \;
       =
       \lim_{k \to \infty} \frac{-b - \log \lceil k t + c \rceil}{J(\sqrt{2 k t})} + \frac{J(-2 + \sqrt{2(1 - \frac{1}{\lceil kt + c \rceil}) (k t + c) + 2})}{J(\sqrt{2 k t})}.
   \end{align*}

   By Assumption \ref{assump: J is sup log} we have $\lim\limits_{k \to \infty} \frac{-b - \log \lceil k t + c \rceil}{J(\sqrt{2 k t})} = \lim\limits_{k \to \infty} \frac{- 2 \log \sqrt{\lceil k t + c \rceil}}{J(\sqrt{2 k t})}= 0$.     We also have 
   \begin{equation*}
       \lim_{k \to \infty} \sqrt{2(1 - \frac{1}{\lceil kt + c \rceil}) (k t + c) + 2} = \sqrt{2 k t + 2c + 2},
   \end{equation*}
   Hence, for large enough $k$, there are fixed $c_2, c_3 \in \mathbb{R}$ such that
   \begin{equation*}
       -2 + \sqrt{2 k t + c_2} \leq -2 + \sqrt{2(1 - \frac{1}{\lceil kt + c \rceil}) (k t + c) + 2} \leq -2 + \sqrt{2 k t + c_3}, \quad \forall k > K.
   \end{equation*}
   Given $J(t)$ is an increasing function, and $\lim_{k \to \infty} \frac{J(-2 + \sqrt{2kt + c_i})}{J(\sqrt{2 k t})} = 1, \; i = 2, 3$ by Lemma \ref{lem: J asymptotics does not change with log} we obtain
   \begin{equation}
       \lim_{k \to \infty}\frac{J(-2 + \sqrt{2(1 - \frac{1}{\lceil kt + c \rceil}) (k t + c) + 2})}{J(\sqrt{2 k t})} = 1.
   \end{equation}
\end{proof}

{
\bibliographystyle{plainnat}
\bibliography{references}
}

\end{document}